\documentclass{amsart}
\usepackage{graphicx}
\title[Blaschke Products as Compositions]{Finite Blaschke Products as Compositions of Other Finite
Blaschke Products}
\author{Carl C. Cowen}%
\address{IUPUI (Indiana University -- Purdue University, Indianapolis), Indianapolis, Indiana 46202-3216}%
\email{ccowen@math.iupui.edu}%
\date{Fall 1974; 16 July 2012}%

   \newcommand{\D}{\ensuremath{\mathbb{D}}}
  \newfont{\caps}{cmcsc10}  
  \newfont{\jour}{cmti10}  
  \newcommand{\vypp}[4]{ {\bf #1}(#2), #3--#4.}
\newcommand{\advmath}[4]{{\jour Advances~in~Math.}\vypp{#1}{#2}{#3}{#4}}
 \newcommand{\pams}[4]{{\jour Proc.~Amer.~Math.~Soc.}\vypp{#1}{#2}{#3}{#4}}
\newcommand{\tams}[4]{{\jour Trans.~Amer.~Math.~Soc.}\vypp{#1}{#2}{#3}{#4}}
  \newcommand{\jfa}[4]{{\jour J. Functional Analysis}\vypp{#1}{#2}{#3}{#4}}
 \newcommand{\jlms}[4]{{\jour J.~London~Math.~Soc.}\vypp{#1}{#2}{#3}{#4}}
  \newcommand{\iu}[4]{{\jour Indiana Univ.~Math.~J.}\vypp{#1}{#2}{#3}{#4}}
\newcommand{\monthly}[4]{{\jour Amer.~Math.~Monthly}\vypp{#1}{#2}{#3}{#4}}

\begin{document}

\begin{abstract}
These notes answer the question ``When can a finite Blaschke product $B$ be
written as a composition of two finite Blaschke products $B_1$ and $B_2$, that is,
$B=B_1\circ B_2$, in a non-trivial way, that is, where the order of each is greater 
than $1$."   It is shown that a group can be computed from $B$ and its local inverses, 
and that compositional factorizations
correspond to normal subgroups of this group.   This manuscript was written in 1974 
but not published because it was pointed out to the author that this was primarily a 
reconstruction of work of Ritt from 1922 and 1923, who reported on work on polynomials.  
It is being made public now because of recent interest in this subject by several mathematicians
interested in different aspects of the problem and interested in applying these ideas to 
complex analysis and operator theory.  
\end{abstract}

\maketitle

\section{Introduction}\label{Intro}

From the point of view of these notes, for a positive integer $n$, a Blaschke product of order $n$
(or $n$-fold Blaschke product) is an $n$-to-one analytic map of the open unit disk, \D, onto itself.
It is well known that such maps are rational functions of order $n$, so have continuous extensions
to the closed unit disk and the Riemann sphere that are $n$-to-one maps of these sets onto themselves,
and have the form 
\[ B(z)=\lambda\prod_{k=1}^n \frac{z-\alpha_j}{1-\overline{\alpha_j}z} \]
for $\alpha_1$, $\alpha_2$, $\cdots$, $\alpha_n$ points of the unit disk and $|\lambda|=1$.  

These notes were my first formal mathematical writing, developed at the beginning of the work on 
my thesis, and were written as a present to my former teacher Professor John Yarnelle on the occasion 
of his retirement from Hanover College, Hanover, Indiana, where I had been a student.   The only
original of these notes was given to Professor Yarnelle (since deceased) in December 1974 and what is
presented here is a scan of the Xerox\texttrademark\ copy I made for myself at that time.  These
notes have never been formally circulated, but they have been shared over the years with several people
and form the basis of the work in my thesis~\cite{Co78a}, especially in Section 2,  my further work 
on commutants of analytic Toeplitz operators then~\cite{Co78b,Co80a,Co80b}, and more recently in work on multiple valued  composition operators~\cite{CG06} and a return to questions of commutants of analytic multiplication operators~\cite{CoW}.  In addition, they have formed the basis of 
my talk ``An Unexpected Group", given to several undergraduate audiences in recent years
starting in 2007 at Wabash College.   In the past few years, more interest has been shown in this 
topic and it seems appropriate to make these ideas public and available to others who are working with 
related topics.  Examples of a revival of interest of Ritt's ideas are in the work of R.~G.~Douglas and D.~Zheng and their collaborators, for example~\cite{DPW, DSZ}, and in the purely function theoretic questions such as the very nice work of Rickards~\cite{Rick} on decomposition of polynomials and 
the paper of Beardon and Ng~\cite{BN}.

The reason this is the first time these notes are being circulated is simple.  In the fall of 1976, I gave a 
talk on this work in the analysis seminar at the University of Illinois at Urbana-Champaign where I was 
most junior of postdocs.  The audience received it politely, and possibly with some interest, so as the 
end of the talk neared, I was feeling good at my first foray into departmental life.  Then, at question time,
the very distinguished Professor Joe Doob asked, ``Didn't Ritt~\cite{Ritt22,Ritt23} do something like this in the 1920's?"   I was devastated and embarrassed and promptly put the manuscript in a drawer, thinking it 
unpublishable.  In retrospect, I probably should have gone to Professor Doob for advice and written it
up for publication with appropriate citations.  Because Ritt's first work was on compositional factorization of 
polynomials, it is somewhat different than this, but it is obvious that the ideas involved apply to polynomials,
Blaschke products, or rational functions more generally.   I believe that many analysts today, as I was 
then, are ignorant of Ritt's work in this arena and, at the very least, his work deserves to be better known.

The remainder of this document is the scan of the original work for Professor Yarnelle from Fall 1974,
a short addendum from a year later that describes the application of these ideas to factorization of an analytic map on the disk into a composition of an analytic function and a finite Blaschke product,
and a short bibliography of some work related to these ideas.

In the original notes, given a (normalized) finite Blaschke product $I$ of order $n$, a group $G_I$
is described as a permutation group of the branches of the local inverses of the Blaschke product $I$
acted on by loops (based at $0$) in a subset of the disk, the disk with $n(n-1)$ points removed.  The main theorem of the notes is the following.\\[-1ex]

\noindent {\bf Theorem 3.1.} \emph{Let $I$ be a finite Blaschke product normalized as above.}

\emph{If ${\mathcal P}$ is a partition of the set of branches of $I^{-1}$ at $0$, $\{g_1, g_2,
\cdots, g_n\}$, that $G_I$ respects, then there are finite Blaschke products $J_{\mathcal P}$ and 
$b_{\mathcal P}$ with the order of $b_{\mathcal P}$ the same as the order of ${\mathcal P}$
so that }
\[ I = J_{\mathcal P} \circ b_{\mathcal P}\] 

\emph{Conversely, if $J$ and $b$ are finite Blaschke products so that $I=J\circ b$, then there is a 
partition ${\mathcal P}_b$ of the set of branches of $I^{-1}$ at $0$ which $G_I$ respects such that
the order of ${\mathcal P}_b$ is the same as the order of $b$.}

\emph{Moreover, if ${\mathcal P}$ and $b$ are as above, then}
\[ {\mathcal P}_{b_{\mathcal P}} = {\mathcal P}_b   \ \ \ \ \mbox{ and } \ \ \ \  b_{{\mathcal P}_b} =  b  \]

It is shown that the compositional factorizations of $G_I$ are associated with normal subgroups of 
$G_I$, but that the association is more complicated than one might hope in that non-trivial normal 
subgroups of $G_I$ can be associated with trivial compositional factorizations of $I$.  However,
the association is strong enough, then if one knows all of the normal subgroups of $G_I$, then one
can construct all possible non-trivial factorizations of $I$ into compositions of finite Blaschke products
and inequivalent factorizations of $I$ as compositions correspond to different normal subgroups of
$G_I$.

The main theorem of the addendum is the following.\\[-1ex]

\noindent {\bf Theorem.} \emph{If $f: \D \mapsto f(\D)$ is analytic and exactly $n$-to-one {\rm [}as a map of the open unit disk onto the image $f(\D)${\rm ]}, then there is a finite Blaschke product $\phi$ and a 
one-to-one function $\widetilde{f}$ so that  $f=\widetilde{f}\circ\phi$.}\\[-1ex]

This result has the obvious corollaries that $f(\D)$ is simply connected and $f'$ has exactly $n-1$
zeros in the disk. 
\newpage
\addtolength{\textheight}{2in}
\addtolength{\topmargin}{-1.5in}
\addtolength{\textwidth}{2.5in}
\addtolength{\oddsidemargin}{-1.5in}
\addtolength{\evensidemargin}{-1.5in}

\includegraphics{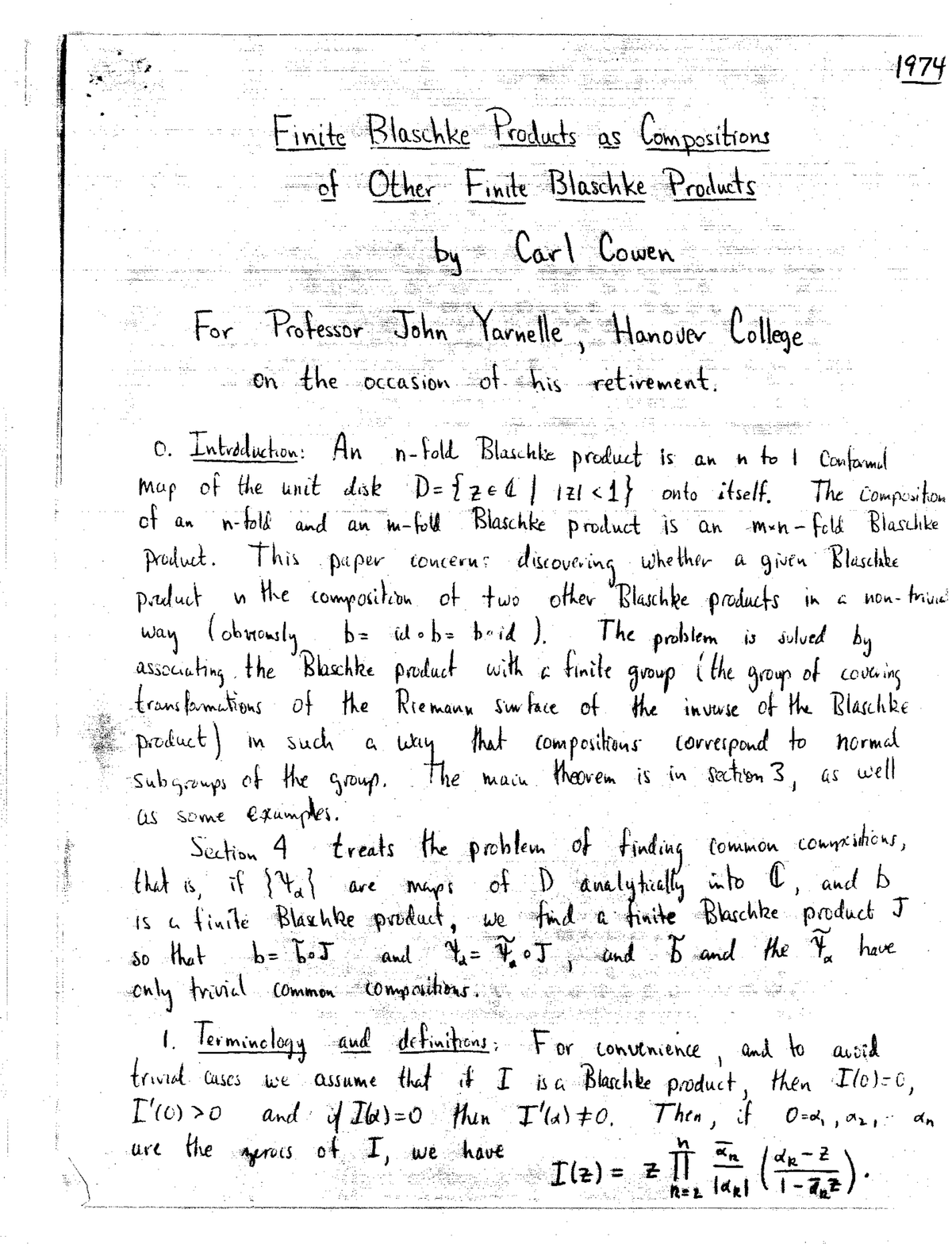}

\newpage

\includegraphics{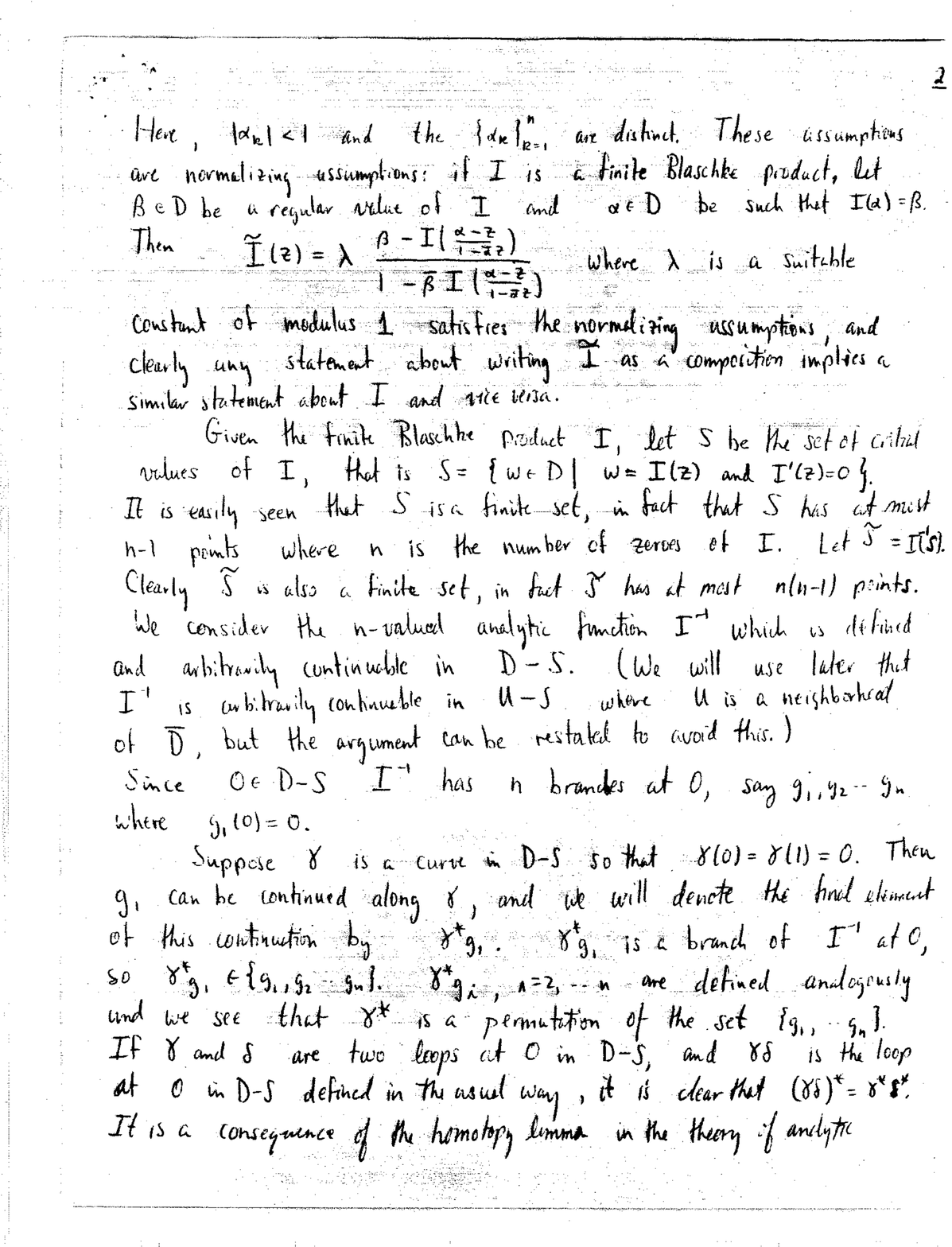}

\newpage

\includegraphics{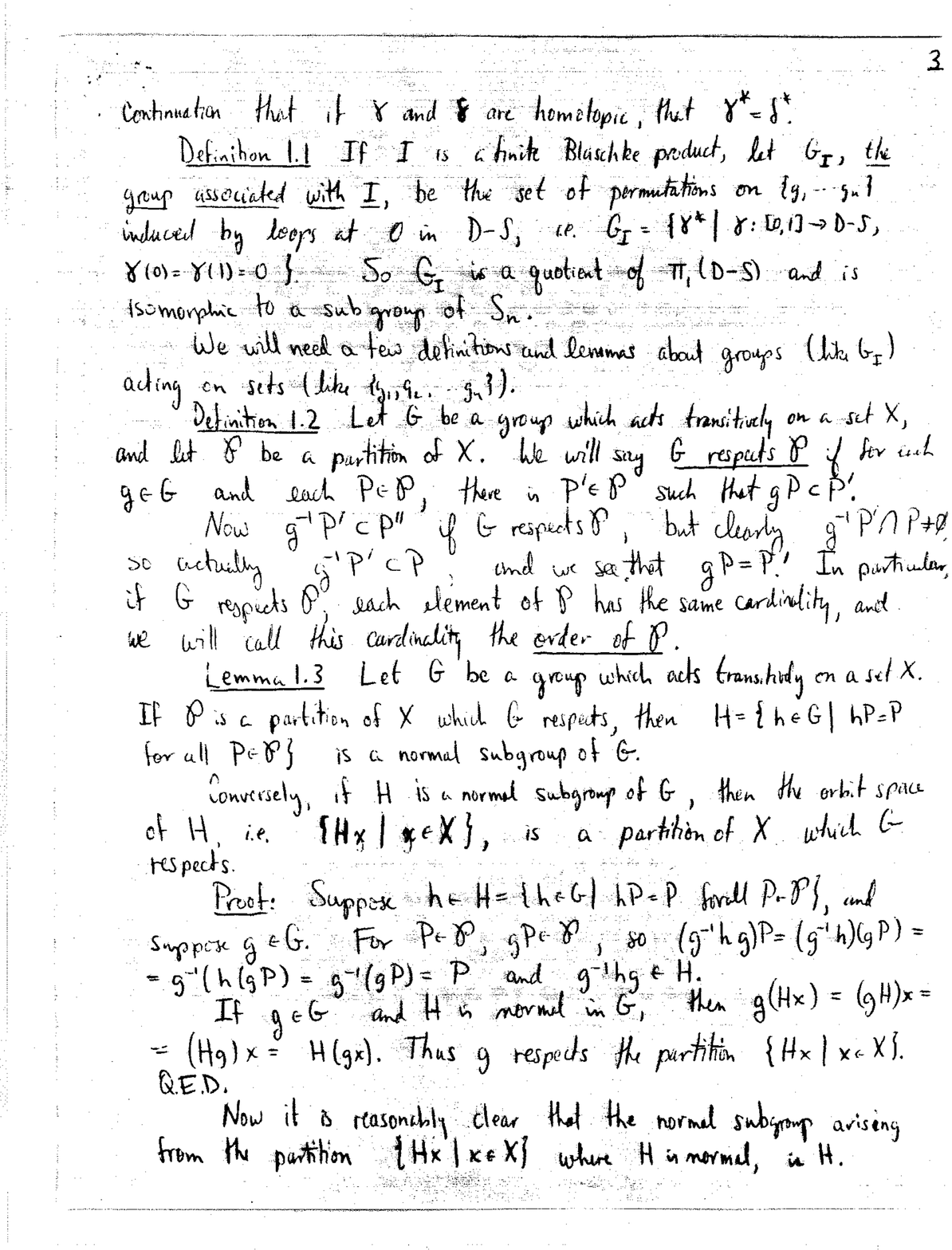}

\newpage

\includegraphics{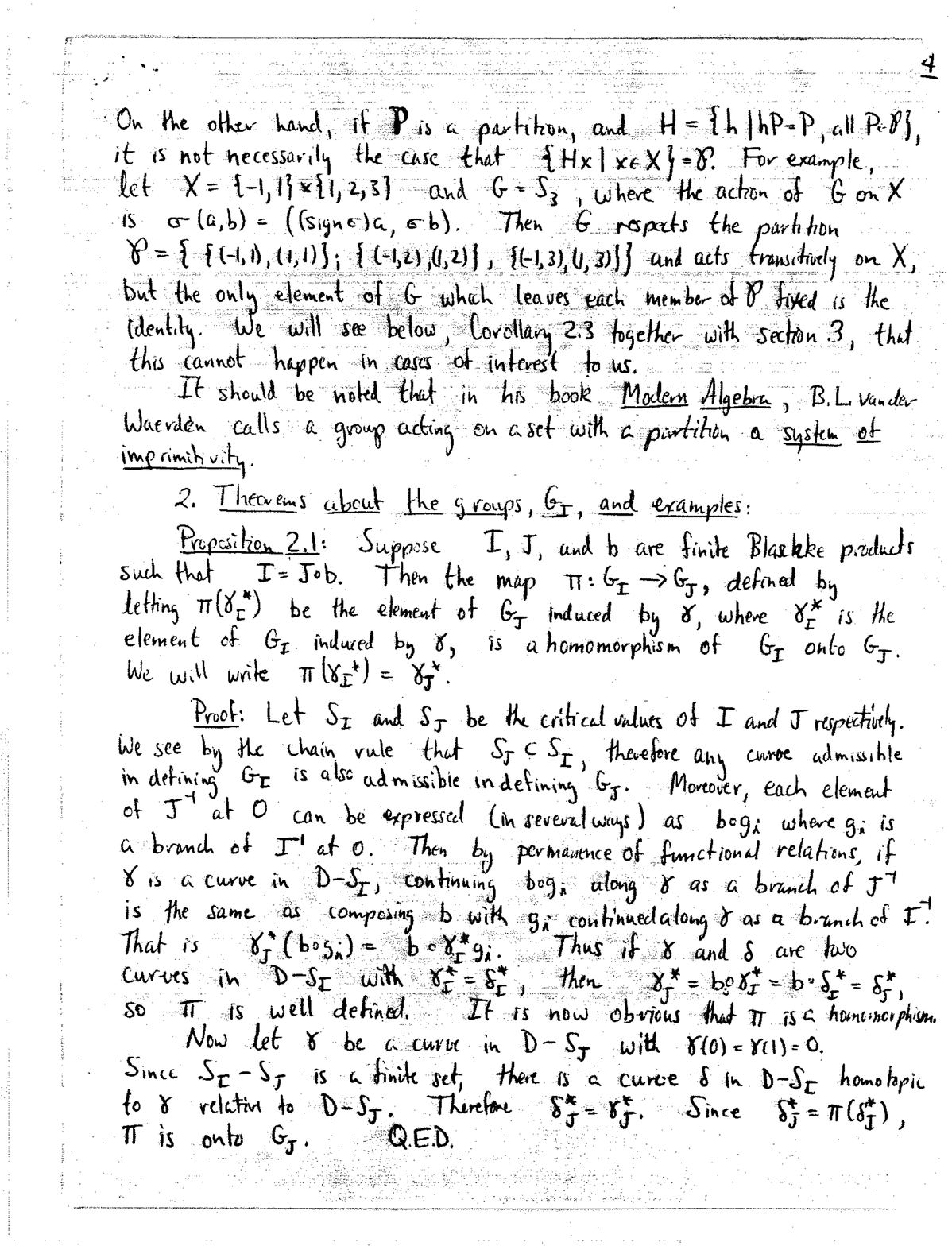}

\newpage

\includegraphics{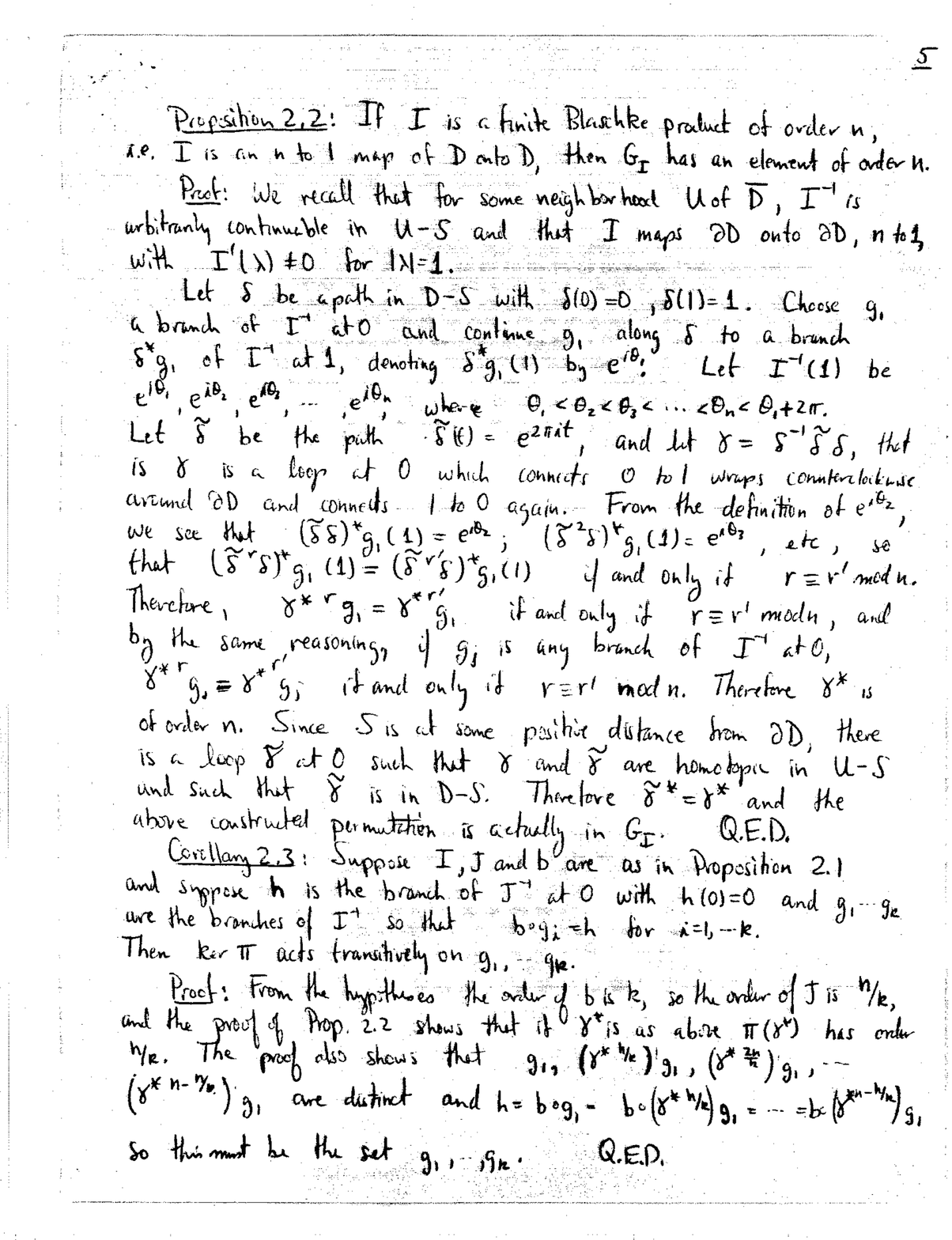}

\newpage

\includegraphics{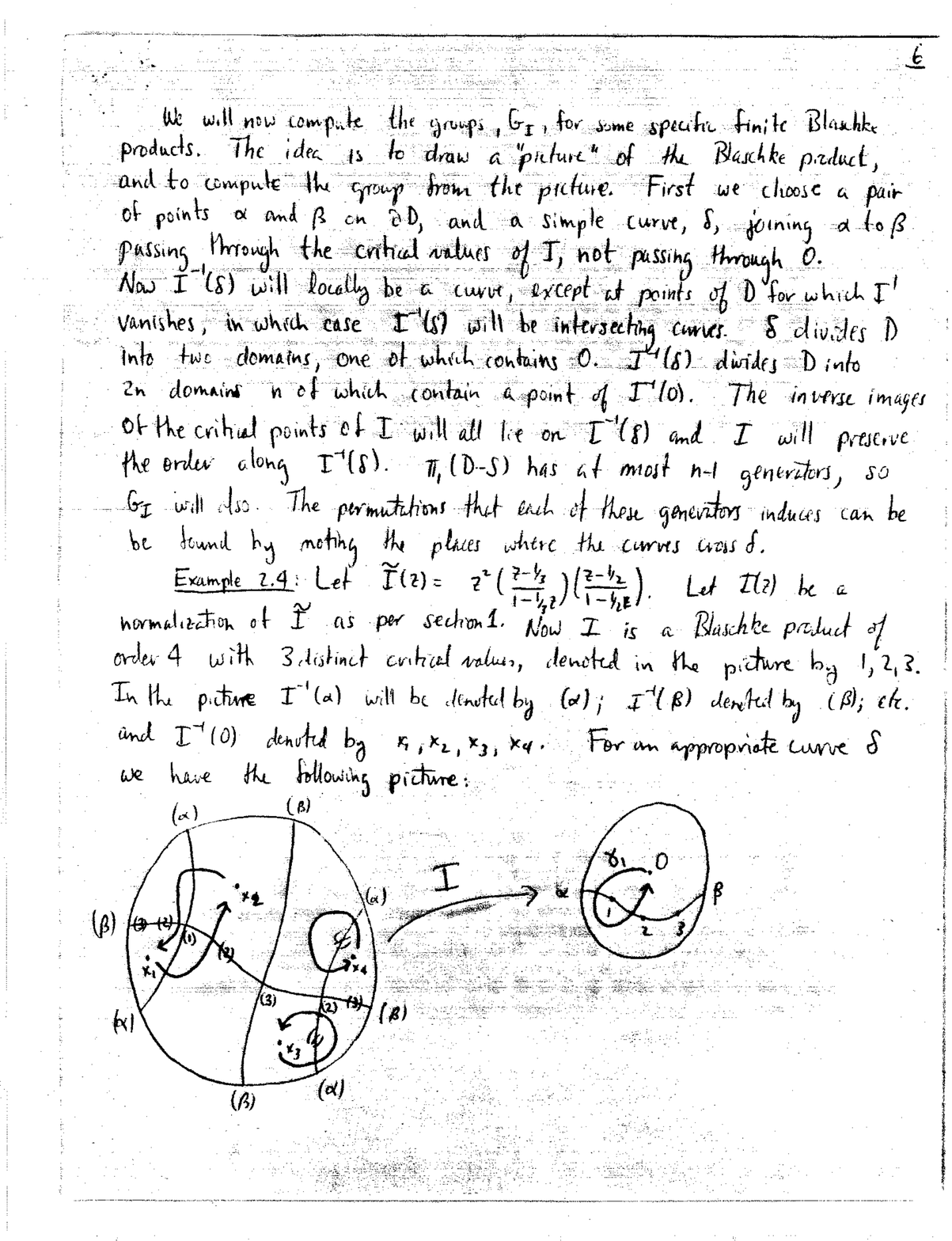}

\newpage

\includegraphics{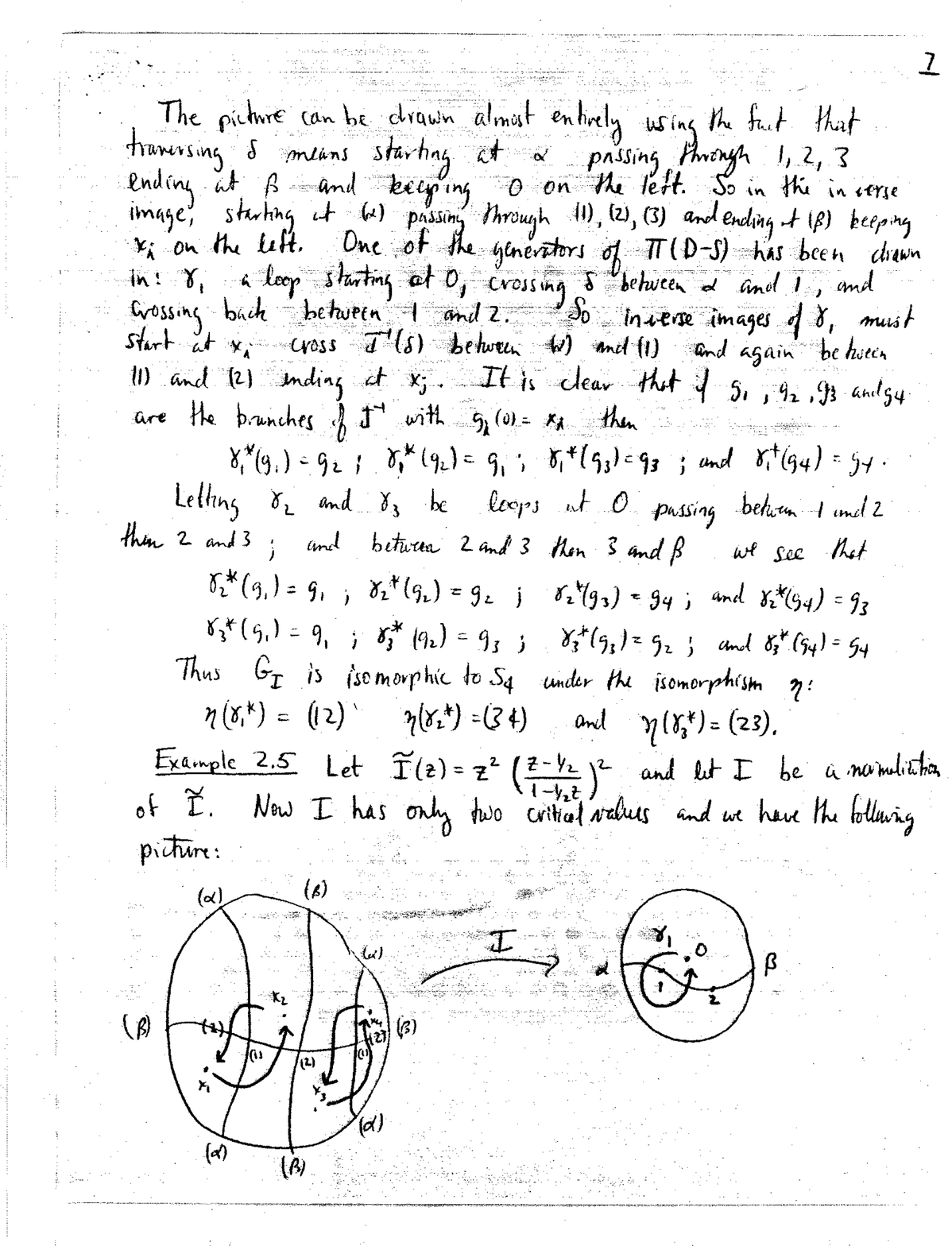}

\newpage

\includegraphics{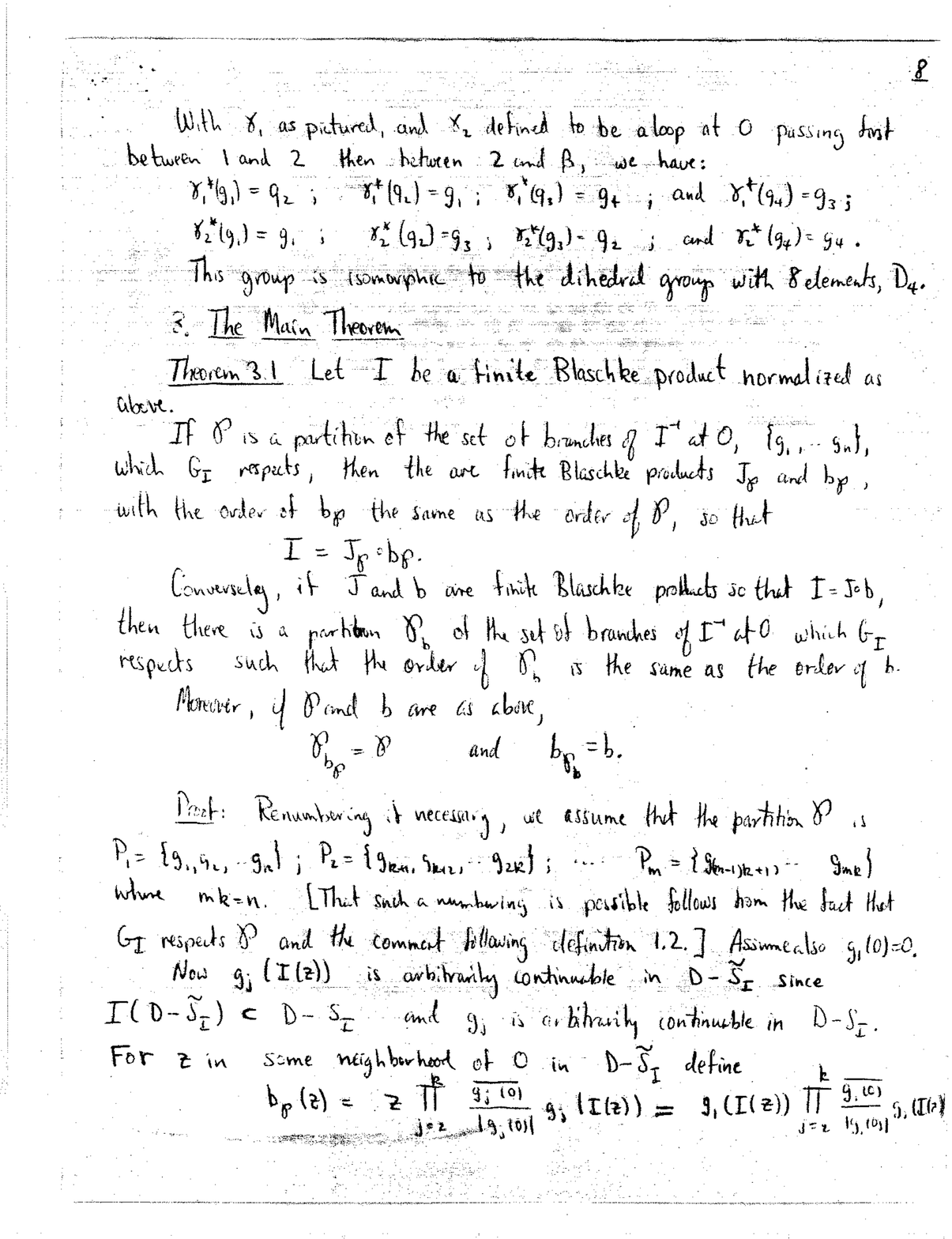}

\newpage

\includegraphics{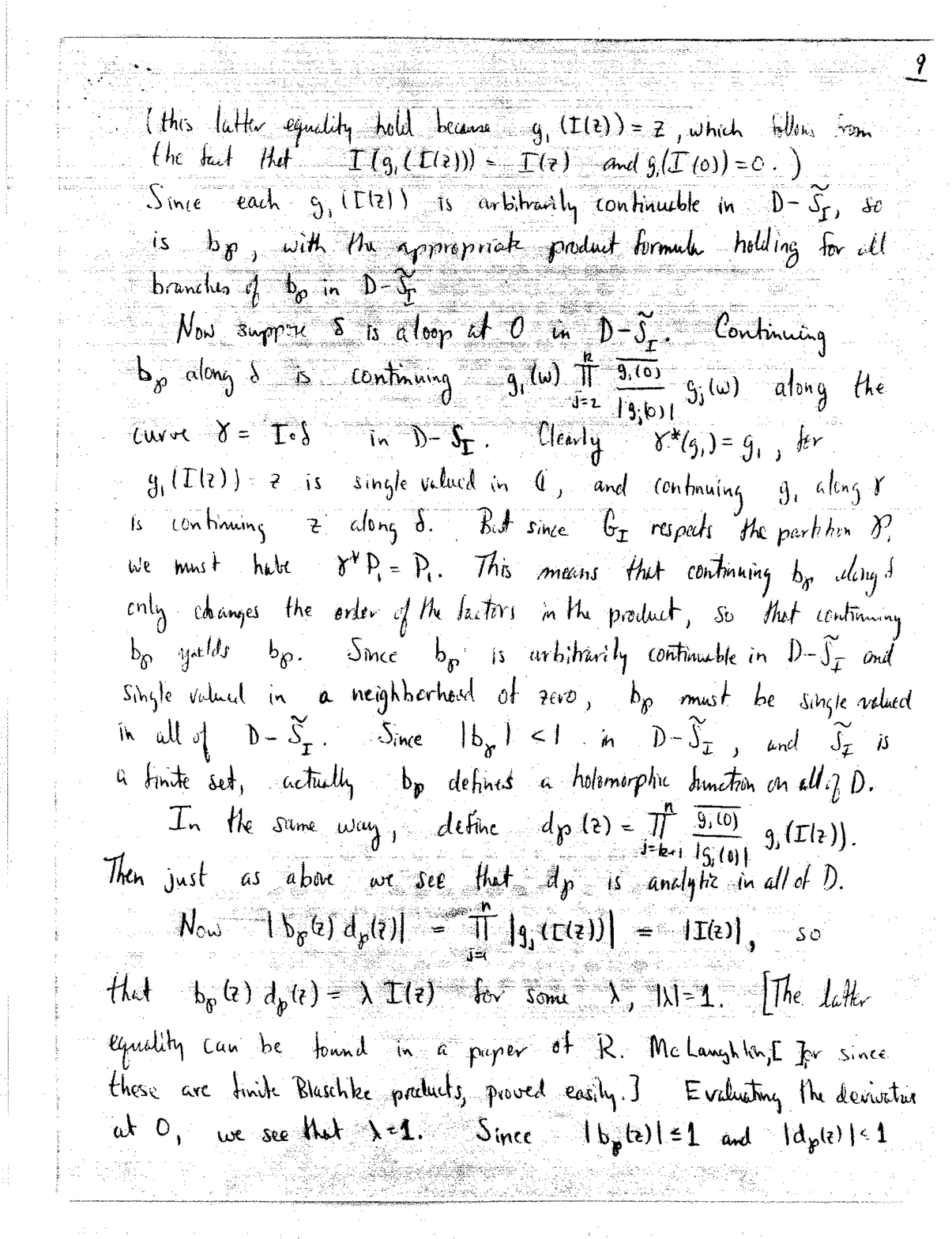}

\newpage

\includegraphics{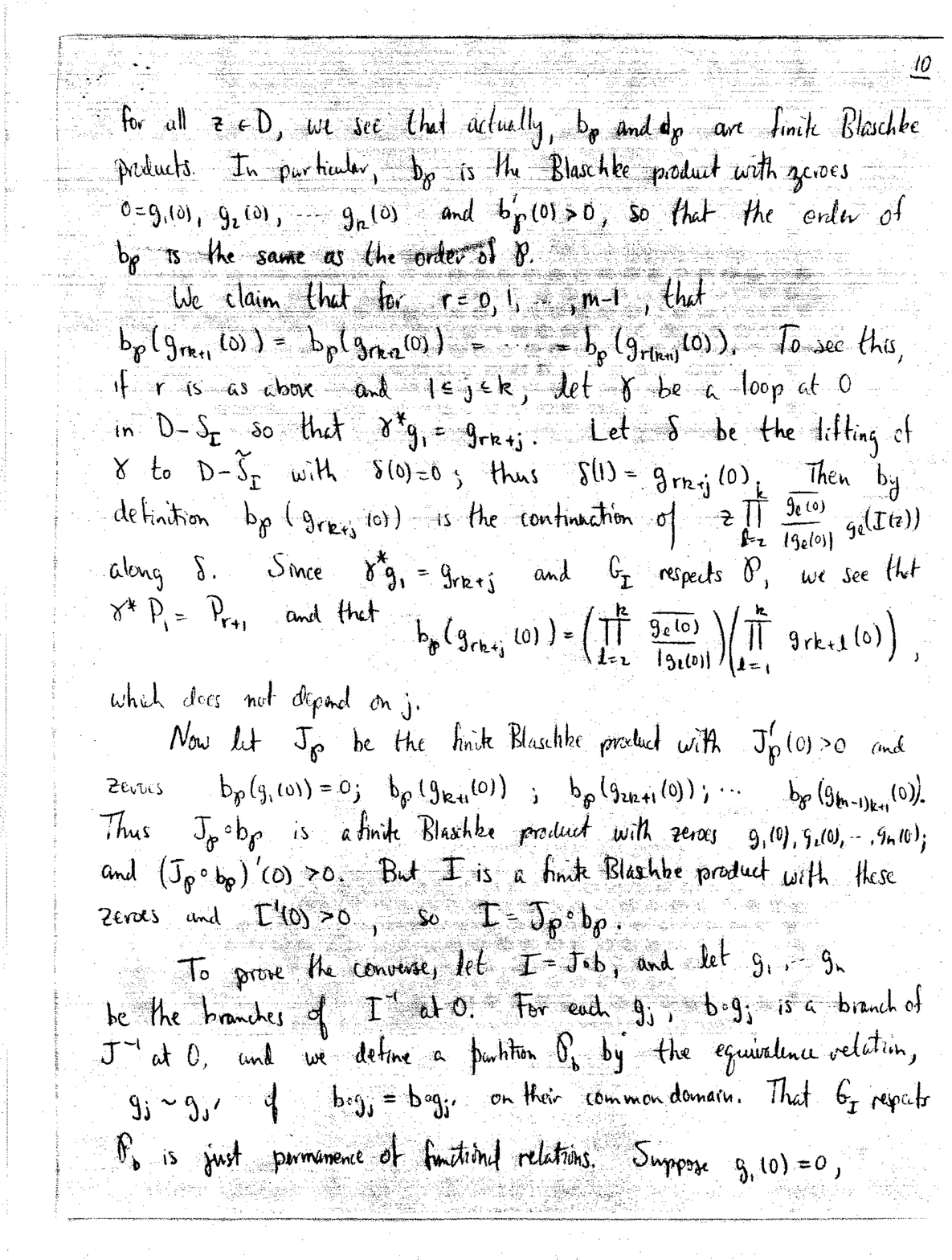}

\newpage

\includegraphics{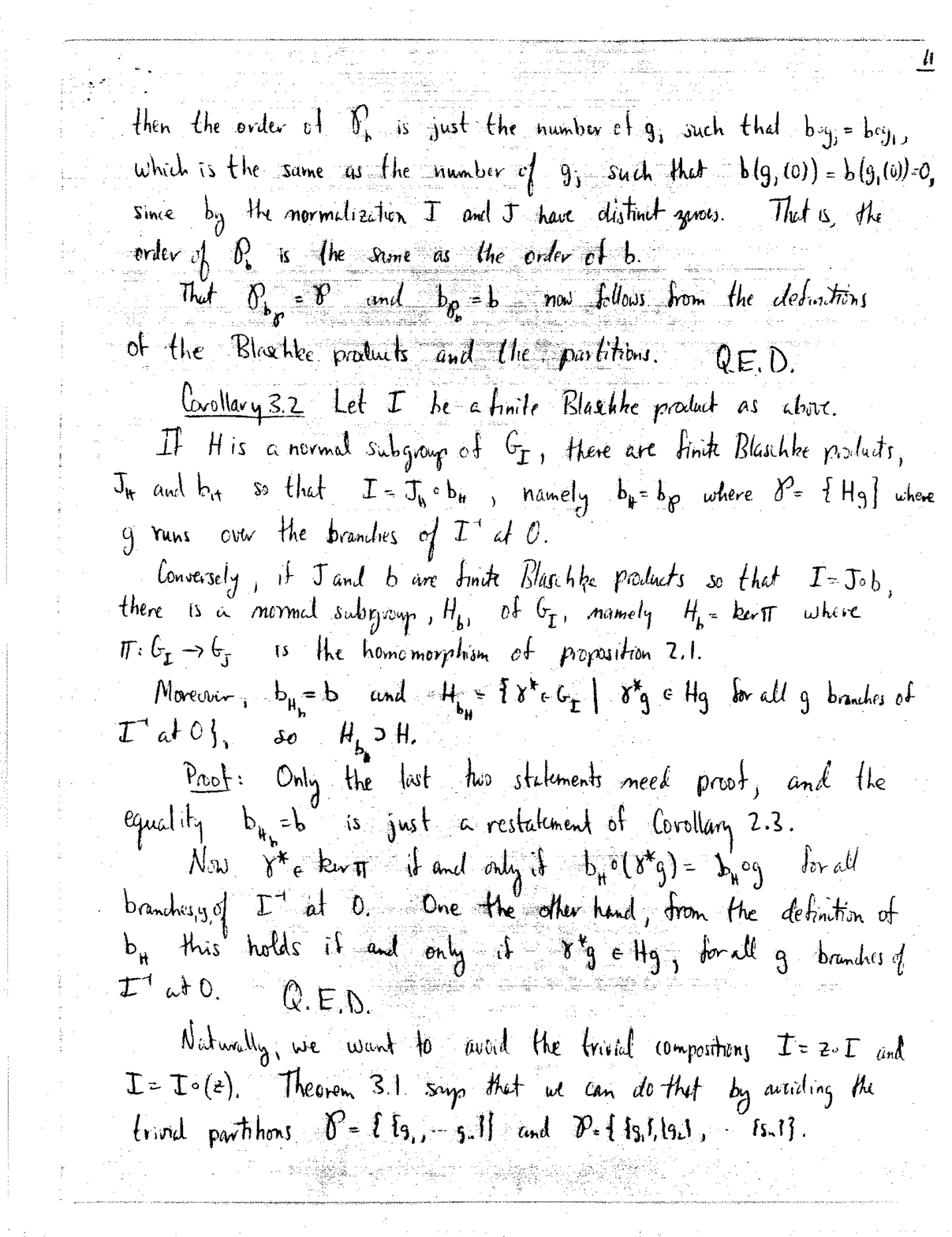}

\newpage

\includegraphics{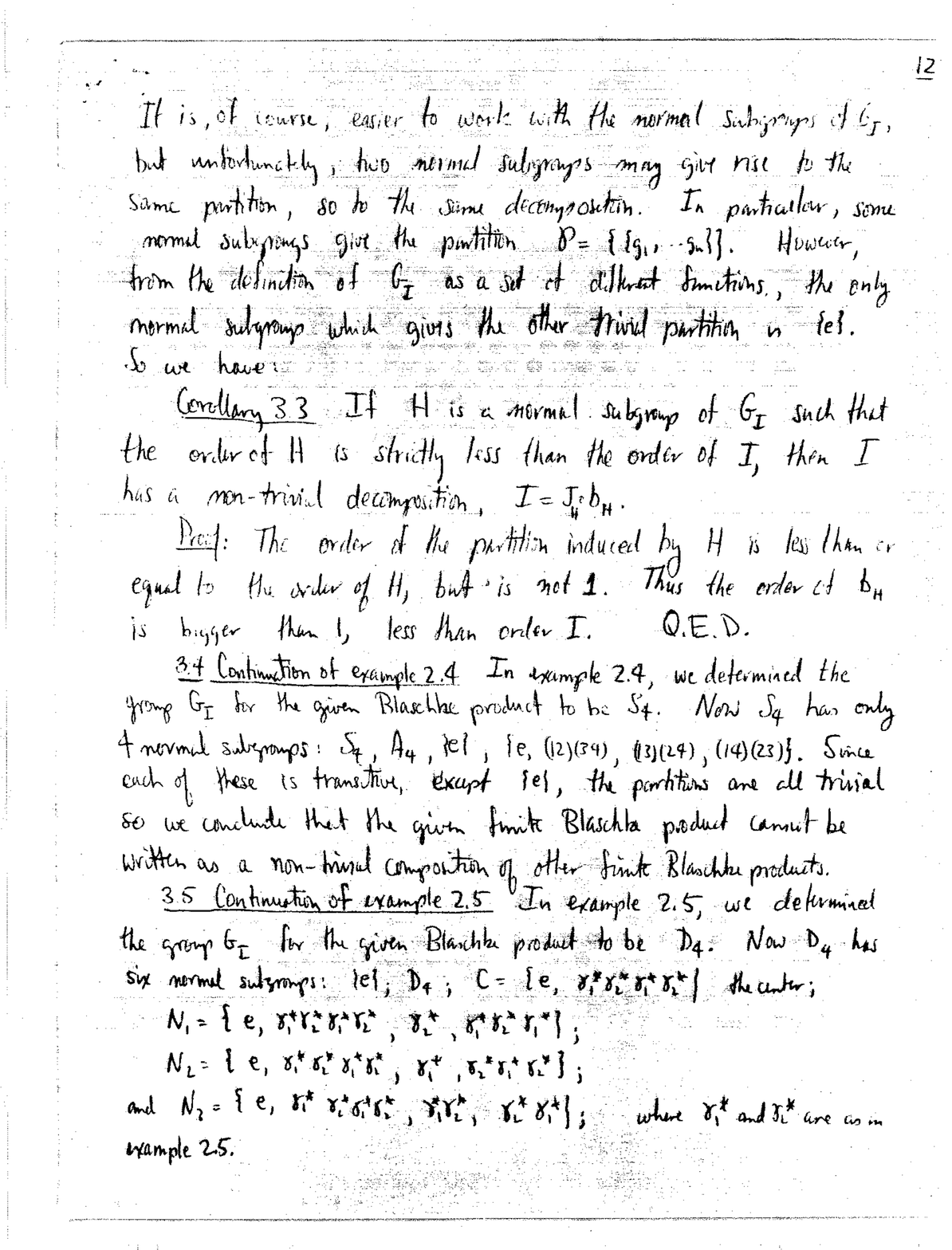}

\newpage

\includegraphics{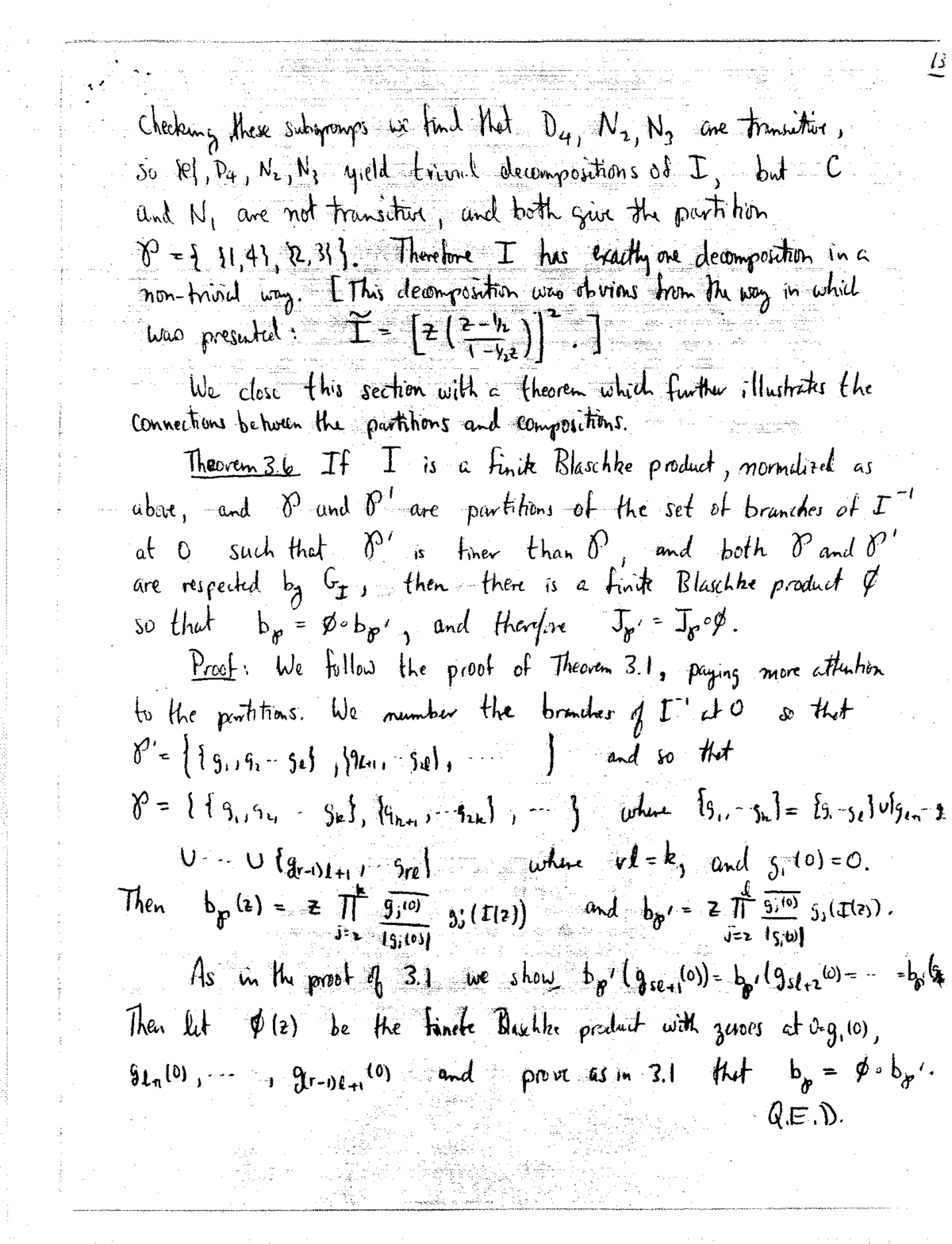}

\newpage

\includegraphics{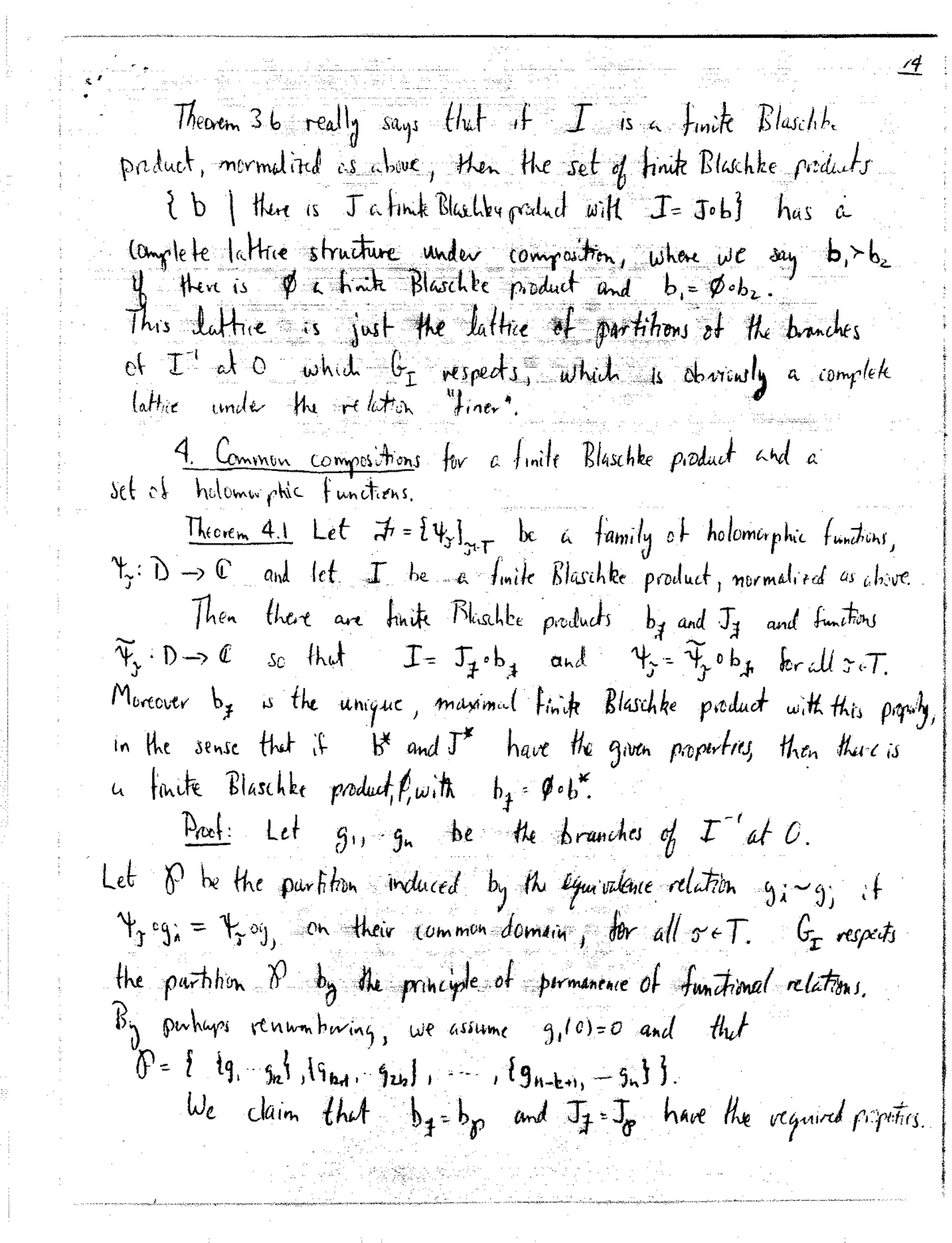}

\newpage

\includegraphics{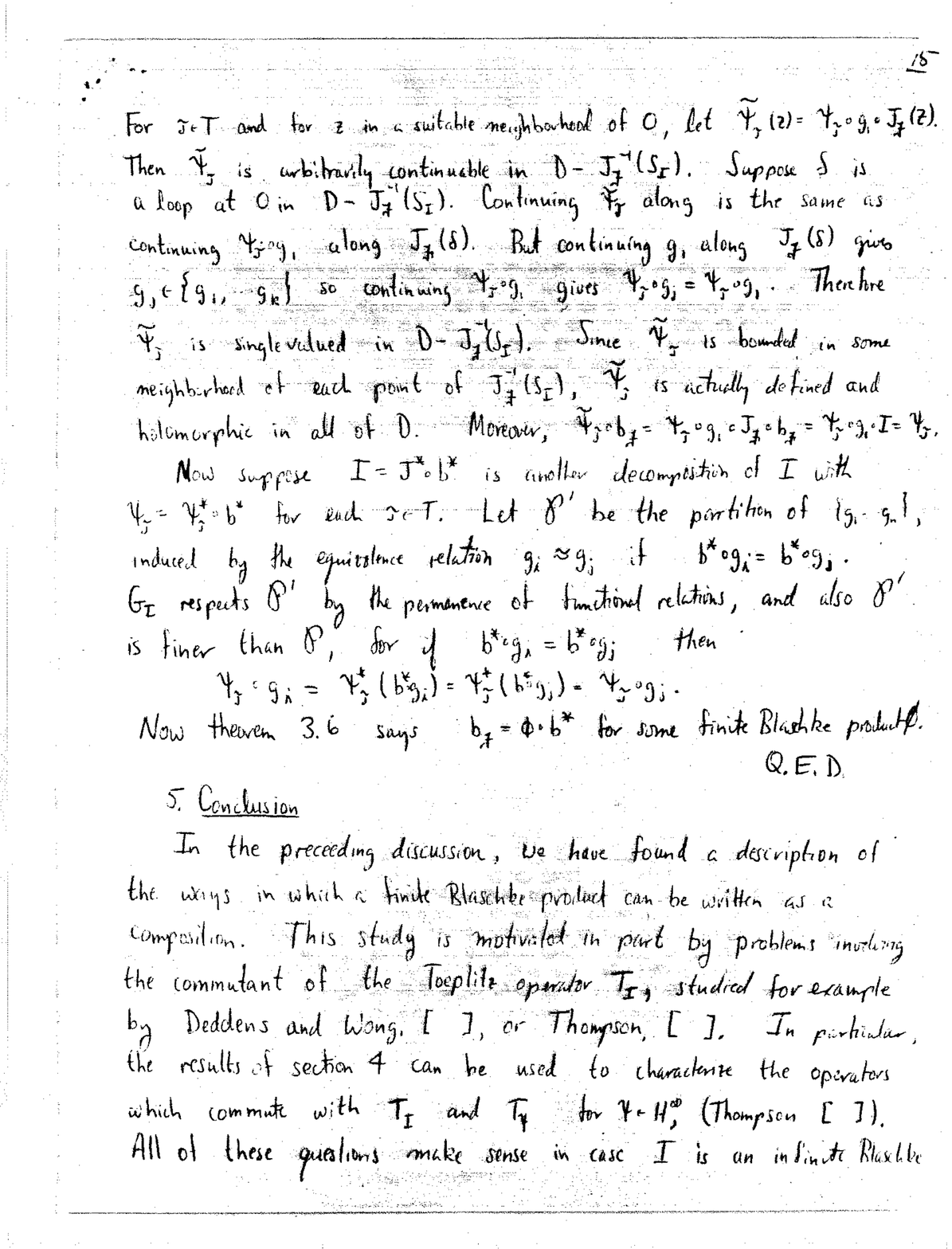}

\newpage

\includegraphics{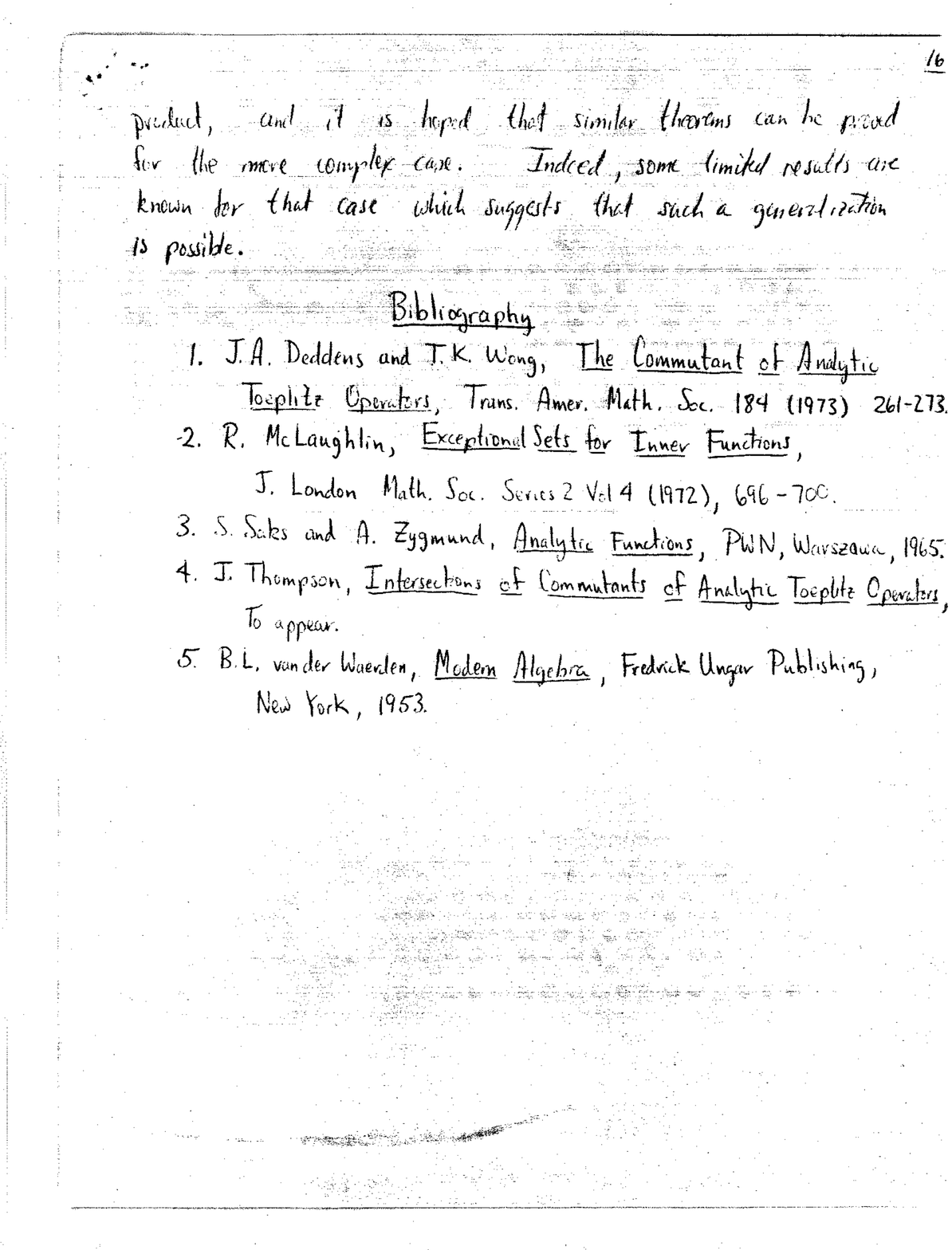}

\newpage
\addtolength{\topmargin}{.2in}

\includegraphics{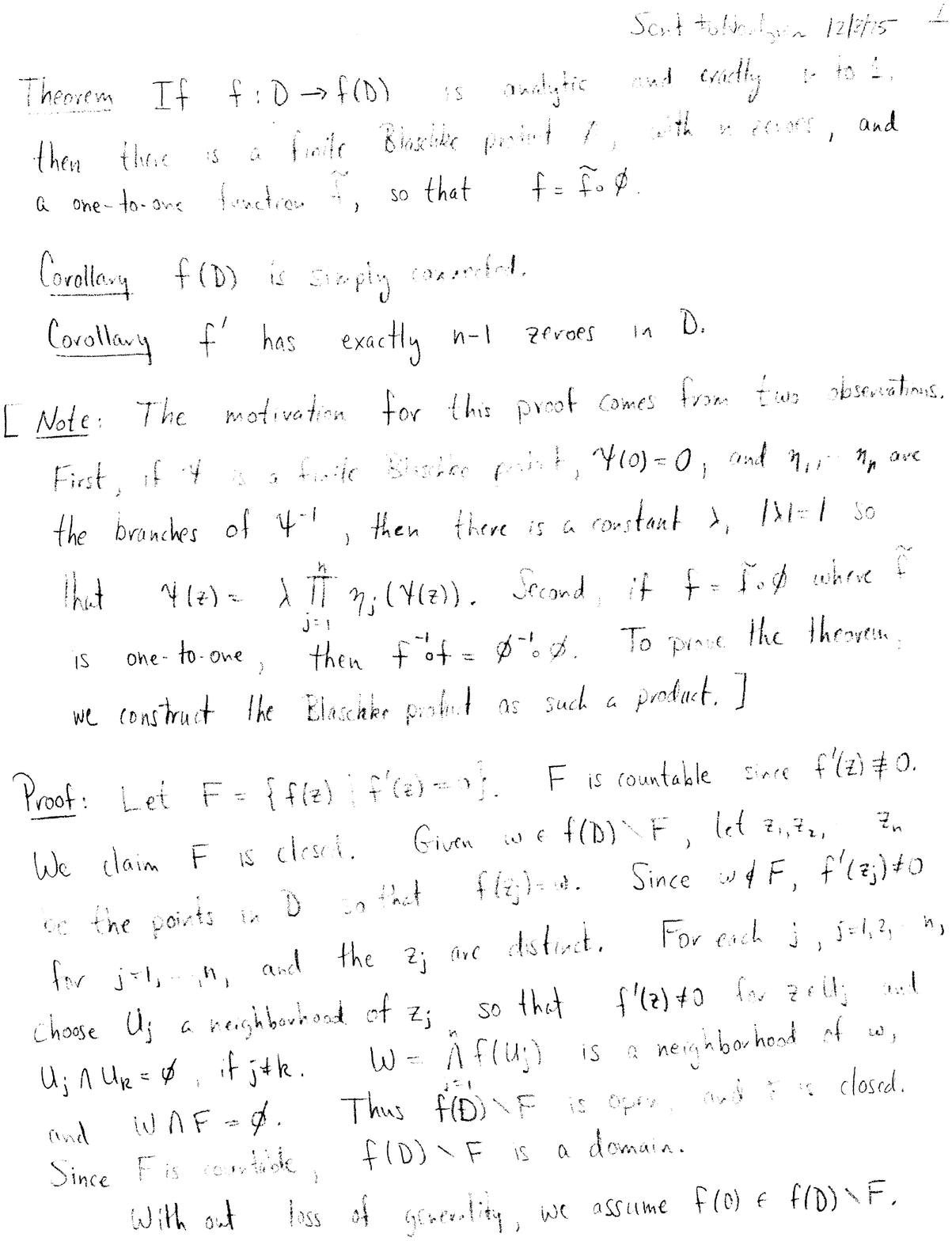}

\newpage

\includegraphics{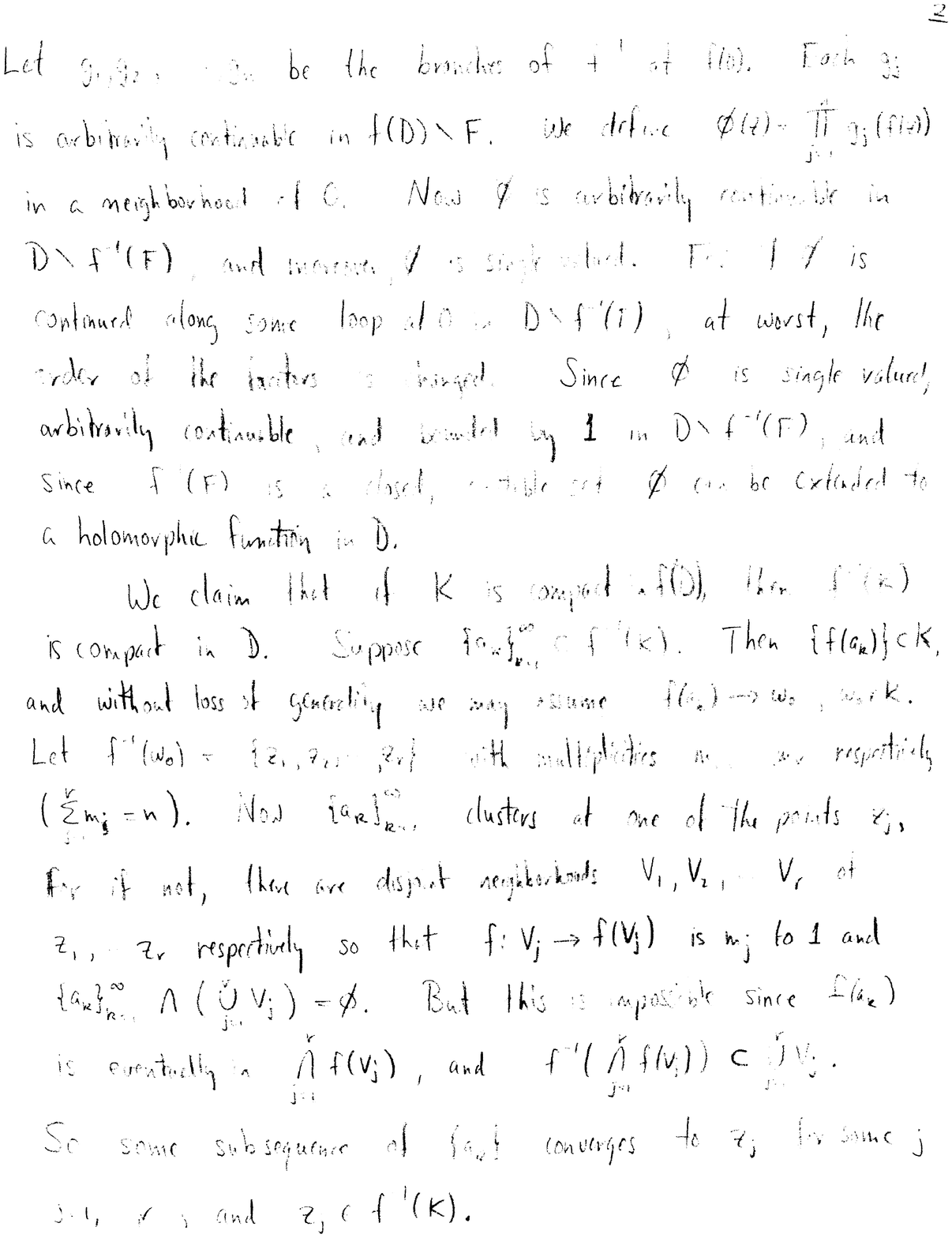}

\newpage

\includegraphics{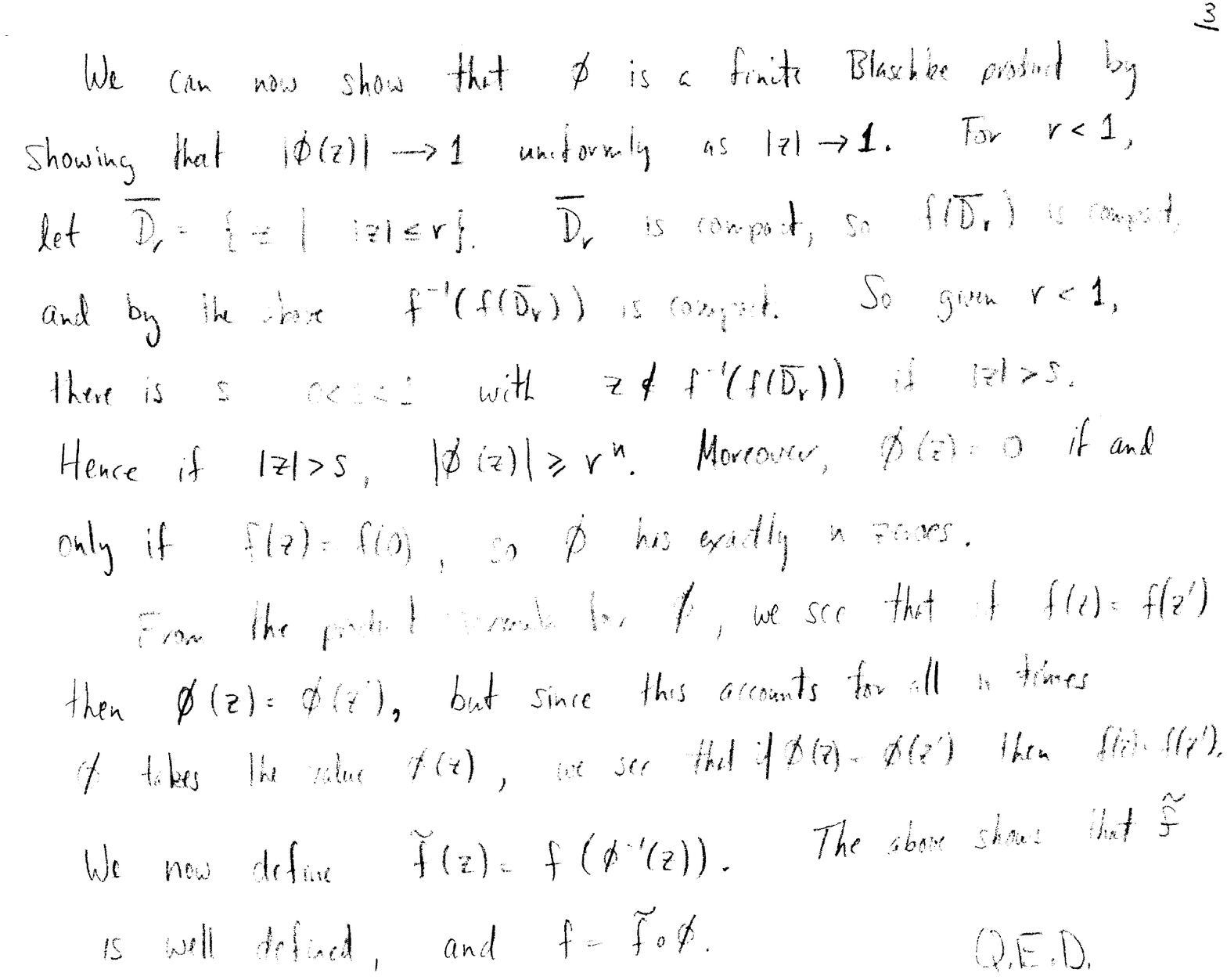}

\newpage

\addtolength{\textheight}{-2in}
\addtolength{\topmargin}{1.3in}
\addtolength{\textwidth}{-2.5in}
\addtolength{\oddsidemargin}{1.5in}
\addtolength{\evensidemargin}{1.5in}


\end{document}